\documentclass[11pt]{article}

\usepackage[a4paper]{geometry}

\usepackage{amsfonts,latexsym,amstext}
\usepackage{amsmath,amscd,amssymb}

\usepackage[latin1]{inputenc}

\usepackage[usenames]{color}

\usepackage{url}

\usepackage{psfrag}

\usepackage{placeins}

\usepackage[dvips]{graphicx}
\usepackage{latexsym,amssymb,amsmath,color,times,bbm}

\newtheorem{thm}{Theorem}[section]
\newtheorem{cor}[thm]{Corollary}
\newtheorem{lem}[thm]{Lemma}
\newtheorem{prop}[thm]{Proposition}
\newtheorem{rem}[thm]{Remark}

\numberwithin{equation}{section}
\newcommand{\norm}[1]{\left\Vert#1\right\Vert}
\newcommand{\abs}[1]{\left\vert#1\right\vert}

\def\cprime{$'$}

\def\sqw{\hbox{\rlap{\leavevmode\raise.3ex\hbox{$\sqcap$}}$%
\sqcup$}}
\def\sqb{\hbox{\hskip5pt\vrule width4pt height6pt depth1.5pt%
\hskip1pt}}

\def\qed{\ifmmode\hbox{\hfill\sqb}\else{\ifhmode\unskip\fi%
\nobreak\hfil
\penalty50\hskip1em\null\nobreak\hfil\sqb
\parfillskip=0pt\finalhyphendemerits=0\endgraf}\fi}
\def\cqfd{\ifmmode\sqw\else{\ifhmode\unskip\fi\nobreak\hfil
\penalty50\hskip1em\null\nobreak\hfil\sqw
\parfillskip=0pt\finalhyphendemerits=0\endgraf}\fi}

\def \proof{{\noindent \bf Proof. }}
\def \eproof{\hbox{ }\hfill$\Box$}

\begin{document}

\renewcommand{\labelitemi}{$\bullet$}
\bibliographystyle{plain}
\pagestyle{headings}
\title{On the uniqueness of solutions to quadratic BSDEs with non-convex generators}
\author{ 
          Philippe Briand
         \\\small Univ. Grenoble Alpes, Univ. Savoie Mont Blanc, CNRS, LAMA, 73000 Chambéry, France
         \\\small \sf philippe.briand@univ-smb.fr
		\and
             Adrien Richou
             \\\small  \small Univ. Bordeaux, IMB, UMR 5251, F-33400 Talence, France.
              \\\small  \sf adrien.richou@math.univ-bordeaux.fr 
              }


\maketitle

\begin{abstract} 
In this paper we prove some uniqueness results for quadratic backward stochastic differential equations without any convexity assumptions on the generator. The bounded case is revisited while some new results are obtained in the unbounded case when the terminal condition and the generator depend on the path of a forward stochastic differential equation. Some of these results are based on strong estimates on $Z$ that are interesting on their own and could be applied in other situations. 
\end{abstract}



\section{Introduction}

In this paper, we consider the following quadratic backward stochastic differential equation (BSDE in short for the remaining of the paper)
\begin{equation}
\label{EDSR}
 Y_t=\xi { +} \int_t^T f(s,Y_s,Z_s)ds -\int_t^T Z_s dW_s, \quad 0 \leqslant t \leqslant T,
\end{equation}
where the generator $f$ has a quadratic growth with respect to $z$. In \cite{Kobylanski-00} Kobylanski studied the case where $\xi$ and the random part of $f$ are bounded. She proved the existence of a solution $(Y,Z)$ such that $Y$ is bounded and she get that this solution is unique amongst solutions $(\tilde{Y},\tilde{Z})$ such that $\tilde{Y}$ is bounded. The unbounded case was { investigated} in \cite{Briand-Hu-06} where authors obtained an existence result. The problem of uniqueness in the unbounded framework was tackled in \cite{Briand-Hu-08,Delbaen-Hu-Richou-09,Delbaen-Hu-Richou-13} by assuming that $f$ is a convex function with respect to $z$. The case of a non-convex generator $f$ was treated in \cite{Richou-12} but uniqueness results where obtained in some classes involving bounds on $Z$.

The main contribution of this paper is to strengthen these uniqueness results. Concerning the bounded case, we are able to expand the class of uniqueness: the bounded solution obtained by { Kobylanski} is unique amongst solutions $(\tilde{Y},\tilde{Z})$ such that $\tilde{Y}$ has a specific exponential moment. In the unbounded framework, we are able to relax the convexity assumption on the generator by assuming that the terminal condition and the random part of the generator depend on the path of a forward stochastic differential equation 
$$X_t = x+\int_0^t b(X_s) ds + \int_0^t \sigma(s,X_s) dW_s.$$
Moreover, the class of uniqueness only involves the process $Y$. To get into the details, two different situations are { investigated}. 
\begin{itemize}
\item When $\sigma$ only { depends} on $s$, we can deal with a terminal condition and a generator that are locally Lipschitz functions of the path of $X$. This uniqueness result relies on a strong estimate on $Z$ given by
 $$|Z_t| \leqslant C(1+\sup_{s \in [0,t]} |X_s|^r), \quad d\mathbb{P} \otimes dt \textrm{ a.e.}$$
This estimate is a { generalization} of an estimate obtained in \cite{Richou-12} in the Markovian framework and is interesting on its own.  
\item When $\sigma$ depends on $X$, we { start by the case of} a terminal condition and a generator that are Lipschitz functions of the path of $X$. In this case, we are able to show that $Z$ is bounded $d\mathbb{P} \otimes dt$ a.e. which is also a new estimate interesting on its own.
\end{itemize}
{
Let us emphasize that, in these two situations, we are able to get a uniqueness result, even if we add a bounded random variable to the terminal condition and a bounded process to the generator. 

\medskip
}

The paper is organized as follows. In section~2, we prove some elementary theoretical uniqueness results that  will be usefull in the following of the article. Finally, Section~3 is devoted to the different frameworks detailed previously: the bounded case and the two different unbouded cases.

Let us close this introduction by giving the notations that we will use in all the paper. For the remaining of the paper, let us fix a nonnegative real number $T>0$. First of all, $(W_t)_{t \in [0,T]}$ is a standard Brownian motion with values in $\mathbb{R}^d$ defined on some complete probability space $(\Omega,\mathcal{F},\mathbb{P})$. $(\mathcal{F}_t)_{t \geqslant 0}$ is the natural filtration of the Brownian motion $W$ augmented by the $\mathbb{P}$-null sets of $\mathcal{F}$. The sigma-field of predictable subsets of $[0,T] \times \Omega$ is denoted by $\mathcal{P}$.

By a solution to the BSDE (\ref{EDSR}) we mean a pair $(Y_t,Z_t)_{t \in [0,T]}$ of predictable processes with values in $\mathbb{R} \times \mathbb{R}^{1\times d}$ such that $\mathbb{P}$-a.s., $t \mapsto Y_t$ is continuous, $t \mapsto Z_t$ belongs to $L^2(0,T)$, $t \mapsto f(t,Y_t,Z_t)$ belongs to $L^1(0,T)$ and $\mathbb{P}$-a.s. $(Y,Z)$ verifies  (\ref{EDSR}). { The terminal condition $\xi$ is $\mathcal{F}_T$-measurable.}

For any real $p\geqslant 1$, $\mathcal{S}^p$ denotes the set of real-valued, adapted and càdlàg processes $(Y_t)_{t \in [0,T]}$ such that
$$\norm{Y}_{\mathcal{S}^p}:=\mathbb{E} \left[\sup_{0 \leqslant t\leqslant T} \abs{Y_t}^p \right]^{1/p} < + \infty.$$
$\mathcal{M}^p$ denotes the set of (equivalent { classes} of) predictable processes $(Z_t)_{t \in [0,T]}$ with values in $\mathbb{R}^{1 \times d}$ such that
$$\norm{Z}_{\mathcal{M}^p}:=\mathbb{E}\left[\left(\int_0^T \abs{Z_s}^2 ds \right)^{p/2}\right]^{1/p} < +\infty.$$
{ We will use the notation $Y^*:=\sup_{0\leqslant t\leqslant T} \abs{Y_t}$ and by $\mathcal{S}^\infty$ we denote the set of adapted càdlàg processes such that $Y^*$ belongs to $L^\infty$.

Let us recall that a continuous local martingale is bounded in mean oscillations if
\begin{equation*}
	||M||_{\text{BMO}_2}= \sup_{\tau} \left|\left|\mathbb{E}[\langle M\rangle_{T} - \langle M \rangle_{\tau}|\mathcal{F}_{\tau}]^{1/2}\right|\right|_{\infty}<\infty
\end{equation*}
where the supremum is taken over all stopping time $\tau \leq T$. We refer to \cite{Kazamaki-94} for further details on BMO-martingales.

Finally, $\mathbb{D}^{1,2}$ stands for the set of random variables $X$ which are differentiable in the Malliavin sense and such that
\begin{equation*}
	\mathbb{E}\left[|X|^2 + \int_0^T |D_s X|^2 \, ds\right] < \infty.
\end{equation*}
Moreover, $\mathbb{L}_{1,2}$ denote the set of real-valued progressively measurable processes $(u_t)_{t \in [0,T]}$ such that
\begin{itemize}
 \item for a.e. $t\in [0,T]$, $u_t \in \mathbb{D}^{1,2}$,
 \item $(t,\omega) \mapsto Du_t(\omega) \in L^2([0,T])$ admits a progressively measurable version,
 \item $\mathbb{E}\left[ \left(\int_0^T |u_t|^2 dt\right)^{1/2}+\left( \int_0^T \int_0^T |D_{\theta}u_t|^2 d\theta dt\right)^{1/2} \right] <+\infty.$
\end{itemize}


}




\section{Some elementary uniqueness results}

We are looking for a uniqueness result for the BSDE
\begin{equation}
 \label{BSDE}
 Y_t = \xi +\int_t^T f(s,Y_s,Z_s)ds - \int_t^T Z_s dW_s, {\quad  0\leq t\leq T,}
\end{equation}
where we assume the following assumptions:
\paragraph{(B1)}
$ $ $f: [0,T] \times \Omega\times \mathbb{R} \times \mathbb{R}^{1 \times d} \rightarrow \mathbb{R}$ is a measurable function with respect to $\mathcal{P} \otimes {\mathcal{B}(\mathbb{R})}\otimes \mathcal{B}(\mathbb{R}^{1 \times d})$. There { exist two constants} $K_y>0$ and $K_z>0$ such that, for all $t \in [0,T]$, $y,y' \in \mathbb{R}$, $z,z' \in \mathbb{R}^{1 \times d}$
\begin{enumerate}
 \item
 $|f(t,y,z)-f(t,y',z)| \leqslant K_{y} |y-y'| \quad \text{a.s.},$
 \item  $\tilde{z} \mapsto f(s,y,\tilde{z})$ is $C^1$ and
 $$|\nabla_z f(s,y,z)-\nabla_z f(s,y,z')| \leqslant K_{z} |z-z'| \quad \text{a.s.}$$
\end{enumerate}

\begin{rem}
\label{rem upper bound gene}
 Since we have
 $$
 f(s,0,z)-f(s,0,0) = z\cdot\nabla_z f(s,0,0)+|z|^2\, \int_0^1 \frac{ z\left(\nabla_z f(s,0,uz)-\nabla_z f(s,0,0) \right)}{|z|^2}\mathbbm{1}_{z \neq 0} du,
 $$
 we can remark that assumption (B1) implies the following upper bound: for all $\eta>0$, for all $s \in [0,T]$, $y \in \mathbb{R}$, $z \in \mathbb{R}^{1 \times d}$, we have
 $$
 |f(s,y,z)| \leqslant |f(s,0,0)| +{ \frac{|\nabla_z f(s,0,0)|^2}{4\eta} }+ K_y |y| + \left(\frac{K_z}{2} + \eta\right)|z|^2 \quad \text{a.s.}$$
\end{rem}

\begin{thm}
\label{resultat unicite 1}
{
Let $p>1$ and $\varepsilon>0$ and let us assume the existence of a solution $(Y,Z)$ to \eqref{BSDE} such that
}
 \begin{equation}
\label{integrabilite expo Z}
 \mathcal{E}_T := e^{\int_0^T \nabla_z f(s,Y_s,Z_s)dW_s-\frac{1}{2} \int_0^T |\nabla_z f(s,Y_s,Z_s)|^2 ds} \in L^p
 \end{equation}
and
  \begin{equation}
\label{integrabilite expo Y}
  \mathbb{E}\left[e^{\frac{2p}{p-1}K_{z}(1+\varepsilon)|Y^*|}\right] <+\infty.
  \end{equation}
Then, this solution is unique amongst solutions  to \eqref{BSDE} such that { the exponential integrability \eqref{integrabilite expo Y} holds true}.
%
%
\end{thm}

\proof
Let us consider $(\tilde{Y},\tilde{Z})$ a solution of \eqref{BSDE} such that
  \begin{equation*}
  \mathbb{E}\left[e^{\frac{2p}{1-p}K_{z}(1+\varepsilon)|\tilde{Y}^*|}\right] <+\infty
  \end{equation*}
and let us denote $\delta Y := \tilde{Y} - Y$ and $\delta Z := \tilde{Z}- Z$. We get
$$
\delta Y_t = 0 + {\int_t^T \left[f(s,\tilde{Y}_s,\tilde{Z}_s) - f(s,{Y}_s,Z_s) \right] ds}-\int_t^T \delta Z_s dW_s
$$
and we can write
$$f(s,\tilde{Y}_s,\tilde{Z}_s) - f(s,{Y}_s,Z_s) = b_s\delta Y_s+\delta Z_s \nabla_z f(s,Y_s,Z_s) + a_s|\delta Z_s|^2$$
with 
$$
b_s := 
          \frac{f(s,\tilde{Y}_s,\tilde{Z}_s)-f(s,Y_s,\tilde{Z}_s)}{\delta Y_s} \, \mathbbm{1}_{|\delta Y_s| >0}
$$
and 
$$
a_s := \int_0^1 \frac{ \delta Z_s \left( \nabla_z f(s,Y_s,Z_s+u\delta Z_s)-\nabla_z f(s,Y_s,Z_s) \right) }{|\delta Z_s|^2} \, \mathbbm{1}_{|\delta Z_s| >0}\, du.
$$
Thanks to assumptions (B1) we know that $|b_s| \leqslant K_{y}$ and $|a_s| \leqslant \frac{K_{z}}{2}$ for all $s \in [0,T]$. Moreover, since \eqref{integrabilite expo Z} is fulfilled, we are allowed to apply Girsanov's theorem: There exists a new probability $\mathbb{Q}$ under which $W^{\mathbb{Q}}:=(W_t-\int_0^t  \nabla_z f(s,Y_s,Z_s) ds )_{t \in [0,T]}$ is a Brownian motion. 
{ 
Thus, we get
$$
\delta Y_t = 0+\int_t^T \left(b_s\delta Y_s + a_s|\delta Z_s|^2\right) ds - \int_t^T \delta Z_s dW^{\mathbb{Q}}_s, \quad 0\leq t\leq T.
$$
For any stopping time $\sigma \leq T$, setting 
\begin{equation*}
	B_s = e^{\int_0^s b_u \mathbbm{1}_{u\geq \sigma} du},
\end{equation*}
we have, from Itô's formula, for any real number $r$,   
\begin{equation*}
	d e^{rB_s \delta Y_s} = r e^{rB_s \delta Y_s} B_s  \left( - b_s \mathbbm{1}_{s < \sigma} \delta Y_s\, ds + \delta Z_s\, dW^{\mathbb{Q}}_s \right) + e^{rB_s \delta Y_s} B_s |\delta Z_s|^2 \left(\frac{r^2}{2}B_s - ra_s\right) ds.
\end{equation*}
In particular, if $\tau \geq \sigma$, since $ra_s \leq |r| K_z/2$,
\begin{align}
	\nonumber
	e^{r \delta Y_\sigma} & = e^{r B_\tau\delta Y_\tau} + \int_\sigma^\tau e^{rB_s \delta Y_s} B_s |\delta Z_s|^2 \left(ra_s - \frac{r^2}{2}B_s \right) ds - \int_\sigma^\tau r e^{rB_s \delta Y_s} B_s \delta Z_s\, dW^{\mathbb{Q}}_s, \\ \label{eq:utile}
	& \leq e^{r B_\tau\delta Y_\tau} + \frac{|r|}{2} \, \int_\sigma^\tau e^{rB_s \delta Y_s} B_s |\delta Z_s|^2 \left(K_z - |r| B_s \right) ds - r\, \int_\sigma^\tau e^{rB_s \delta Y_s} B_s \delta Z_s\, dW^{\mathbb{Q}}_s.
\end{align}
}
For the remaining of the proof we set $\eta=\left((4K_y)^{-1} \log (1+\varepsilon)\right)\wedge T$ which implies in particular that { $e^{-K_y \eta} \geqslant (1+\varepsilon)^{-1/4}$}. For any $n \in \mathbb{N}^*$ we define the stopping time 
{
$$
\tau_n := \inf \left\{ t \in [T-\eta,T] \left| \int_{T-\eta}^t | e^{K_{z}\sqrt{1+\varepsilon}e^{\int_{T-\eta}^s b_u du} \delta Y_s}e^{\int_{t \wedge \tau_n}^{s} b_u du} \delta Z_s|^2 ds>n\right. \right\}.
$$
Je propose de "simplifier" un peu la définition de $\tau_n$ car ce qui est devant est continu
\begin{equation*}
	\tau_n := \inf \left\{ t \in [T-\eta,T] \;\left|\;  \int_{T-\eta}^t |\delta Z_s|^2 ds>n\right. \right\}.
\end{equation*}

Let $t\in[T-\eta,T]$ and let us use the inequality \eqref{eq:utile} with $\sigma = t\wedge\tau_n$, $\tau = \tau_n$ and $r = K_z (1+\varepsilon)^{1/2}$. For $\sigma\leq s \leq \tau$, 
\begin{equation*}
	(1+\varepsilon)^{1/4}\geqslant e^{K_y\eta}\geqslant B_s \geqslant e^{-K_y\eta} \geqslant (1+\varepsilon)^{-1/4}.
\end{equation*}
Thus $|r| B_s \geqslant K_z(1+\varepsilon)^{1/4} \geqslant K_z$ and \eqref{eq:utile} gives
\begin{equation}
\label{inegalite expo delta Y}
e^{K_{z}\sqrt{1+\varepsilon}\,\delta Y_{t \wedge \tau_n}} \leqslant \mathbb{E}_t^{\mathbb{Q}} [e^{K_{z}\sqrt{1+\varepsilon} B_{\tau_n} \delta Y_{\tau_n}}]\leqslant \mathbb{E}_t^{\mathbb{Q}} [e^{K_{z}(1+\varepsilon)^{3/4} |\delta Y_{\tau_n}|}].
\end{equation}
By applying H\"older inequality and by using \eqref{integrabilite expo Y} for $Y$ and $\tilde{Y}$,  we can remark that
\begin{align}
 \nonumber \mathbb{E}^{\mathbb{Q}}\left[e^{ K_{z}(1+\varepsilon) |\delta Y_{\tau_n}| }\right] & = \mathbb{E}\left[ \mathcal{E}_Te^{K_{z}(1+\varepsilon) |\delta Y_{\tau_n} |}\right]\\
 \nonumber &\leqslant \mathbb{E} \left[\mathcal{E}_T^p\right]^{1/p} \mathbb{E} \left[e^{\frac{p}{p-1}K_{z}(1+\varepsilon) |\delta Y_{\tau_n}|}\right]^{\frac{p-1}{p}}\\
 \label{inegalite holder} &\leqslant \mathbb{E} \left[\mathcal{E}_T^p\right]^{1/p} \mathbb{E} \left[e^{\frac{2p}{p-1}K_{z}(1+\varepsilon) |Y^*|}\right]^{\frac{p-1}{2p}}\mathbb{E} \left[e^{\frac{2p}{p-1}K_{z}(1+\varepsilon) |\tilde{Y}^*|}\right]^{\frac{p-1}{2p}}<+\infty.
\end{align}
Thus, $(e^{ K_{z}(1+\varepsilon)^{3/4} |\delta Y_{\tau_n}|})_{n \in \mathbb{N}}$ is uniformly integrable under $\mathbb{Q}$. Since we clearly have that $\tau_n \rightarrow T$ a.s. and $\delta Y_{\tau_n} \rightarrow 0$ a.s. when $n \rightarrow +\infty$, we get 
$$
\mathbb{E}_t^{\mathbb{Q}} [e^{K_{z}(1+\varepsilon)^{3/4} |\delta Y_{\tau_n}| }] \rightarrow 1 \quad \textrm{a.s.}
$$
By taking $n \rightarrow +\infty$ in \eqref{inegalite expo delta Y} we finally obtain that $\tilde{Y}_t \leqslant Y_t$ a.s. for all $t \in [T-\eta,T]$. By the same argument (the quadratic term in \eqref{eq:utile} depends on $|r|$), we can also derive the inequality
\begin{equation*}
e^{-K_{z} \sqrt{1+\varepsilon}\delta Y_{t \wedge \tau_n}} \leqslant \mathbb{E}_t^{\mathbb{Q}} [e^{K_{z}(1+\varepsilon)^{3/4} |\delta Y_{\tau_n}|}], \quad \forall t \in [T-\eta,T],
\end{equation*}
which gives us that $\tilde{Y}_t \geqslant Y_t$ a.s. for all $t \in [T-\eta,T]$. Finally, $\mathbb{E}[\sup_{s \in [T-\eta,T]}|\delta Y_s|^2]=0$ since $Y$ and $\tilde{Y}$ are continuous a.s. It is clear that we can iterate the proof on intervals $[T-(k+1)\eta,T-k\eta]\cap [0,T]$ for $k \in\mathbb{N}^*$ to get that $\mathbb{E}[\sup_{s \in [0,T]}|\delta Y_s|^2]=0$,  As usual it is sufficient to apply It\^o formula to $\delta Y$ to obtain that $\mathbb{E}\left[\int_0^T |\delta Z_s|^2 ds\right] = 0$ which concludes the proof.
\eproof

}

By using same arguments we can also obtain two other versions of this result.
\begin{thm}
\label{resultat unicite 2}
We assume the existence of a solution $(Y,Z)$ to \eqref{BSDE}  such that
 \begin{equation}
\label{integrabilite expo Z 2}
 \mathcal{E}_T := e^{\int_0^T \nabla_z f(s,Y_s,Z_s)dW_s-\frac{1}{2} \int_0^T |\nabla_z f(Z_s)|^2 ds} \in \bigcap_{p>1} L^p
 \end{equation}
and
  \begin{equation}
\label{integrabilite expo Y 2}
  e^{K_{z} |Y^*|} \in \bigcap_{p>1} L^p.
  \end{equation}
Then, this solution is unique amongst solutions $(Y,Z)$ to \eqref{BSDE} such that 
\begin{equation}
\label{integrabilite expo Y 2 classe solution}
e^{K_{z} |Y^*|} \in \bigcup_{p>1} L^p.
\end{equation}
\end{thm}
\begin{thm}
\label{resultat unicite 3}
We assume the existence of a solution $(Y,Z)$ to \eqref{BSDE}  such that
 \begin{equation}
\label{integrabilite expo Z 3}
 \mathcal{E}_T = e^{\int_0^T \nabla_z f(s,Y_s,Z_s)dW_s-\frac{1}{2} \int_0^T |\nabla_z f(Z_s)|^2 ds} \in \bigcup_{p>1} L^p
 \end{equation}
and
  \begin{equation}
\label{integrabilite expo Y 3}
  e^{K_{z}|Y^*|} \in \bigcap_{p>1} L^p.
  \end{equation}
Then, this solution is unique amongst solutions $(Y,Z)$ to \eqref{BSDE} { for which}
\begin{equation}
\label{integrabilite expo Y 3 classe solution}
e^{K_{z} |Y^*|} \in \bigcap_{p>1} L^p.
\end{equation}
\end{thm}
\proof
The proof of Theorem \ref{resultat unicite 2} and Theorem \ref{resultat unicite 3} are overall similar to the previous one. We only sketch the proof of Theorem \ref{resultat unicite 2}, the proof of Theorem \ref{resultat unicite 3} following same lines: we consider $(\tilde{Y},\tilde{Z})$ a solution of \eqref{BSDE} such that
  \begin{equation}
   \label{integrabilite expo Ytilde 2}
  e^{K_{z} |\tilde{Y}^*|} \in \bigcup_{p>1} L^p,
  \end{equation}
  and we show that $Y=\tilde{Y}$ a.s.
  The only difference is in the inequality \eqref{inegalite holder}: instead of applying Cauchy-Schwarz inequality, we use H\"older inequality to get, for any $r>1$ and $p>1$,
\begin{align*}
  \mathbb{E}^{\mathbb{Q}}\left[e^{ K_{z} (1+\varepsilon) \delta Y_{\tau_n}}\right] & \leqslant \mathbb{E} \left[\mathcal{E}_T^p\right]^{1/p} \mathbb{E} \left[e^{\frac{rp}{(r-1)(p-1)}K_{z} (1+\varepsilon) |Y^*|}\right]^{\frac{(r-1)(p-1)}{rp}}\mathbb{E} \left[e^{\frac{rp}{p-1}K_{z} (1+\varepsilon) |\tilde{Y}^*|}\right]^{\frac{p-1}{rp}}.
\end{align*}
Then, by taking $p>1$ large enough, $r>1$ small enough and $\varepsilon>0$ small enough we obtain that 
$$\mathbb{E} \left[\mathcal{E}_T^p\right]^{1/p} \mathbb{E} \left[e^{\frac{rp}{(r-1)(p-1)}K_{z}e^{K_y T} (1+\varepsilon) |Y^*|}\right]^{\frac{(r-1)(p-1)}{rp}}\mathbb{E} \left[e^{\frac{rp}{p-1}K_{z}e^{K_y T} (1+\varepsilon) |\tilde{Y}^*|}\right]^{\frac{p-1}{rp}} <+\infty$$
thanks to \eqref{integrabilite expo Z 2}, \eqref{integrabilite expo Y 2} and \eqref{integrabilite expo Ytilde 2}. The remaining of the proof stays the same.
\eproof

\section{Applications to particular frameworks}

\subsection{The bounded case}
Since the seminal paper of Kobylanski \cite{Kobylanski-00} it is now well known that we have existence and uniqueness of a solution $(Y,Z) \in \mathcal{S}^{\infty} \times \mathcal{M}^2$ to \eqref{BSDE} when $\xi$ and $(f(s,0,0))_{s \in [0,T]}$ are bounded. We are now able to extend the uniqueness to a larger class of solution.

\begin{prop}
\label{prop unicite cas borne}
  We assume that 
  \begin{equation*}
   M:=|\xi|_{L^{\infty}} + \left|\int_0^T |f(s,0,0)|+ \sup_{y \in \mathbb{R}} |\nabla_z f(s,y,0)|ds\right|_{L^{\infty}} <+\infty.
  \end{equation*}
Then there exists $q>1$ that depends only on $M$, $K_y$ and $K_z$ such that the BSDE \eqref{BSDE} admits a unique solution $(Y,Z)$ satisfying 
\begin{equation*}
 \mathbb{E}\left[ e^{2K_zq|Y^*|}\right]<+\infty.
\end{equation*}
In particular, the BSDE \eqref{BSDE} admits a unique solution $(Y,Z)$ satisfying
\begin{equation*}
 e^{K_z |Y^*|} \in \bigcap_{p >1} L^p.
\end{equation*}
\end{prop}
\proof
Thanks to Kobylanski \cite{Kobylanski-00} we know that the BSDE \eqref{BSDE} admits a unique solution $(Y,Z) \in \mathcal{S}^{\infty} \times \mathcal{M}^2$ and this solution satisfies
$$|Y^*|_{L^{\infty}}+\left| \int_0^. Z_s dW_s \right|_{BMO} <+\infty.$$
It implies that 
\begin{eqnarray*}
\left| \int_0^. \nabla_z f (s,Y_s,Z_s) dW_s \right|_{BMO}& =& \left| \sup_{\tau \in [0,T]} \mathbb{E}_{\tau} \left[ \int_{\tau}^T |\nabla_z f (s,Y_s,Z_s)|^2 ds\right]^{1/2}\right|_{L^{\infty}}\\
&\leqslant & \left|\int_0^T\sup_{y \in \mathbb{R}} |\nabla_z f(s,y,0)|ds\right|_{L^{\infty}}+K_z \left| \int_0^. Z_s dW_s \right|_{BMO} <+\infty.
\end{eqnarray*}
Then, the reverse H\"older inequality (see e.g. \cite{Kazamaki-94}) implies that there exists $p^*>1$ such that 
$$\mathcal{E}_T \in \bigcap_{1 \leqslant p < p^*} L^p.$$
Finally we just have to apply Theorem \ref{resultat unicite 1}: for any $\varepsilon>0$ we have the uniqueness of the solution amongst solutions $(Y,Z)$ that satisfy 
\begin{equation}
  \mathbb{E}\left[e^{\frac{2p^*}{p^*-1}K_{z}(1+\varepsilon)|Y^*|} \right] <+\infty.
  \end{equation}
\eproof

\begin{rem}
 It is possible to have an estimate of the exponent $q$ appearing in Proposition \ref{prop unicite cas borne}. Indeed, following the proof, the exponent $q$ is a function of $p^*$ and, using the proof of Theorem 3.1 in \cite{Kazamaki-94}, this exponent $p^*$ is given by
 $$p^* := \phi^{-1}\left(\left| \int_0^. \nabla_z f (s,Y_s,Z_s) dW_s \right|_{BMO}\right),\quad \textrm{with } \phi : p \mapsto \left( 1+\frac{1}{q^2}\log \frac{2q-1}{2q-2}\right)^{1/2} -1.$$
 Moreover, $\left| \int_0^. \nabla_z f (s,Y_s,Z_s) dW_s \right|_{BMO}$ is bounded by an explicit function of $ \left| \int_0^. Z_s dW_s \right|_{BMO}$ and we have some estimates of this last quantity, see for example \cite{Briand-Elie-13}.
\end{rem}

\subsection{A first unbounded case}

In this subsection we consider an SDE with an additive noise
\begin{equation}
  \label{SDE 1}
  X_t = x + \int_0^t b(X_s)ds+\int_0^t \sigma(s) dW_s, \quad 0 \leqslant t \leqslant T,
\end{equation}
where $b$ and $\sigma$ satisfy classical assumptions:
\paragraph{(F1)}
\begin{enumerate}
 \item $b: \mathbb{R}^d \rightarrow \mathbb{R}^d$ is a Lipschitz function: for all $(x,x') \in \mathbb{R}^d \times \mathbb{R}^d$ we have $|b(x)-b(x')| \leqslant K_b |x-x'|$.
 \item $\sigma:[0,T] \rightarrow \mathbb{R}^{d \times d}$ is a bounded measurable function.
\end{enumerate}

We want to study the following BSDE
\begin{equation}
 \label{BSDE path dep 1}
 Y_t = \xi + h((X_s)_{s \in [0,T]})+\int_t^T f(s,Y_s,Z_s)+g((X_{u \wedge s})_{u \in [0,T]},Y_s,Z_s)ds-\int_t^T Z_s dW_s,
\end{equation}
with $h: \mathcal{C}([0,T],\mathbb{R}^d) \rightarrow \mathbb{R}$, $f: [0,T] \times \Omega\times \mathbb{R} \times \mathbb{R}^{1 \times d} \rightarrow \mathbb{R}$ and $g: \mathcal{C}([0,T],\mathbb{R}^d) \times \mathbb{R} \times \mathbb{R}^{1 \times d} \rightarrow \mathbb{R}$ some measurable functions with respect to $\mathcal{B}(\mathcal{C}([0,T],\mathbb{R}^d))$,  $\mathcal{P} \otimes \mathcal{B}(\mathbb{R}^d) \otimes \mathcal{B}(\mathbb{R}^{1 \times d})$ and $\mathcal{B}(\mathcal{C}([0,T],\mathbb{R}^d))\otimes\mathcal{B}(\mathbb{R}^d) \otimes \mathcal{B}(\mathbb{R}^{1 \times d})$ .
We will assume following assumptions:
\paragraph{(B2)}
\begin{enumerate}
 \item $$|\xi|_{L^{\infty}}+\left|\int_0^T |f(s,0,0)|ds\right|_{L^{\infty}}+ \left| \sup_{s \in[0,T],y \in \mathbb{R}, \mathbf{x} \in C([0,T],\mathbb{R}^d)} |\nabla_z f (s,y,0)|+ |\nabla_z g(\mathbf{x},y,0)|  \right|_{L^{\infty}} <+\infty$$ 
 \item there exists $C>0$ such that, for all $t\in [0,T]$, $y \in \mathbb{R}$ and $z \in \mathbb{R}^{1\times d}$,
 $$|f(t,y,z)| \leqslant C.$$
 \item There exist $K_h>0$, $K_g>0$ and $ r \in [0,1)$ such that, for all $\mathbf{x}, \mathbf{\tilde{x}} \in C([0,T],\mathbb{R}^d)$, $y \in \mathbb{R}$, $z \in \mathbb{R}^{1 \times d}$,
 $$|h(\mathbf{x}) -h({\mathbf{\tilde{x}}})| \leqslant K_h(1+|\mathbf{x} |_{\infty}^r+|\mathbf{\tilde{x}}|_{\infty}^r )|\mathbf{x} - \mathbf{\tilde{x}}|_{\infty},$$
 $$|g(\mathbf{x},y,z) -g({\mathbf{\tilde{x}},y,z})| \leqslant K_g(1+|\mathbf{x} |_{\infty}^r+|\mathbf{\tilde{x}}|_{\infty}^r )|\mathbf{x} - \mathbf{\tilde{x}}|_{\infty},$$
 \item (B1) holds true for $(s,y,z) \mapsto f(s,y,z)$ and $(s,y,z) \mapsto g((X_{u \wedge s})_{u \in [0,T]},y,z)$.
\end{enumerate}

\begin{prop}
\label{prop:majorationZ}
 We consider the path-dependent framework and so we assume that $\xi=0$ and $f=0$. We also assume that Assumptions (F1)-(B2) hold true. Then there exists a solution $(Y,Z)$ of the path-dependent BSDE \eqref{BSDE path dep 1} in $\mathcal{S}^2 \times \mathcal{M}^2$ such that,
 \begin{equation}
 \label{est:Z path dependent bruit additif}
 |Z_t| \leqslant C(1+\sup_{s \in [0,t]} |X_s|^r) \quad \quad d\mathbb{P} \otimes dt \textrm{ a.e.}
 \end{equation}
\end{prop}
\proof
The Markovian case was already treated in \cite{Richou-12}. The idea is to generalize this result to the discrete path dependent case, as in \cite{Hu-Ma-04}, and then pass to the limit to obtain the general path dependent case. Since the only novelty is the gathering of known methods and results, we will only sketch the proof.

1. First of all, we start by localizing the generator $g$ to obtain a Lipschitz continuous generator. Let us consider $\rho_N$ a regularized version of the projection on the centered Euclidean ball of radius $N$ in $\mathbb{R}^{1 \times d}$ such that $|\rho_N| \leqslant N$, $|\nabla \rho_N| \leqslant 1$ and $\rho_N(x)=x$ when $|x| \leqslant N-1$. We denote $(Y^N,Z^N) \in \mathcal{S}^2 \times \mathcal{M}^2$ the unique solution of the BSDE
$$Y_t^N = h((X_s)_{s \in [0,T]}) + \int_t^T g^N((X_{u \wedge s})_{u \in [0,T]},Y^N_s,Z^N_s)ds-\int_t^T Z^N_s dW_s$$
with $g^N = g(.,.,\rho_N(.))$. In the remaining of the proof we will see how to prove that \eqref{est:Z path dependent bruit additif} is satisfied by $Z^N$ with a constant $C$ that does not depend on $N$. Let us remark that this is sufficient to conclude since it is quite standard to show that $(Y^N,Z^N)$ is a Cauchy sequence in $\mathcal{S}^2 \times \mathcal{M}^2$ and that the limit is solution of \eqref{BSDE path dep 1},  by using for example a linearization argument and the uniform estimate on $Z^N$. For the reading convenience we will skip the superscript $^N$ in the following. 

2. We approximate $h$ and the random part of $g$ by some discrete functions: by a mere generalization of \cite{Zhang-04} there exists a family $\Pi=\{\pi\}$ of partitions of $[0,T]$ and some families of discrete functionals $\{h^{\pi}\}$, $\{g^{\pi}\}$ such that, for any $\pi \in \Pi$, assuming $\pi: 0=t_0<...<t_n =T$, we have
\begin{itemize}
 \item $h^{\pi} \in \mathcal{C}^{\infty}_b (\mathbb{R}^{d(n+1)})$ and $g^{\pi}(.,y,z)\in \mathcal{C}^{\infty}_b (\mathbb{R}^{d(n+1)})$ for all $(y,z) \in \mathbb{R} \times \mathbb{R}^{1 \times d}$,
 \item $\sum_{i=0}^n |\partial_{x_i} h^{\pi}(x)| \leqslant K_h(1+ 2\sup_{0 \leqslant i \leqslant n} |x_i|^r)$ for all $x = (x_0,...,x_n) \in \mathbb{R}^{d(n+1)}$,
 \item $\sum_{i=0}^n |\partial_{x_i} g^{\pi}(x,y,z)| \leqslant K_g(1+ 2\sup_{0 \leqslant i \leqslant n} |x_i|^r)$ for all $x = (x_0,...,x_n) \in \mathbb{R}^{d(n+1)}$ and $(y,z) \in \mathbb{R} \times \mathbb{R}^{1 \times d}$,
 \item $\lim_{|\pi| \rightarrow 0} |h^{\pi}(x_{t_0},...,x_{t_n})-h(x)|=0$, for all $x \in \mathcal{C}([0,T],\mathbb{R}^d)$,
 \item $\lim_{|\pi| \rightarrow 0} |g^{\pi}(x_{t_0},...,x_{t_n},y,z)-g(x,y,z)|=0$, for all $x \in \mathcal{C}([0,T],\mathbb{R}^d)$ and $(y,z) \in \mathbb{R} \times \mathbb{R}^{1 \times d}$.
\end{itemize}
Let us emphasize that $K_h$ and $K_g$ do not depend on $N$ and $\pi$. We firstly assume that $g$ is smooth enough with respect to $y,z$ and $b$ is smooth enough with respect to $x$, then we have the representation 
\begin{equation*}
 Z_t^{\pi} = \nabla Y_t^{\pi} \nabla X_t^{-1} \sigma(t) \quad \quad \forall t \in [0,T],
\end{equation*}
where
\begin{align*}
 \nabla Y_t^{\pi} &= \sum_{i=1}^n \nabla^i Y_t^{\pi} \mathbbm{1}_{[t_{i-1},t_i)}(t) + \nabla^n Y_{T^-}^{\pi}  \mathbbm{1}_{\{T\}}(t) \\
 \nabla^i Y_t^{\pi} &= \sum_{j \geqslant i} \partial_{x_j} h^{\pi} \nabla X_{t_j}+\int_t^T \sum_{j \geqslant i} \partial_{x_j}g^{\pi} \nabla X_s + \partial_y g^{\pi} \nabla^i Y_s^{\pi} + \partial_z g^{\pi} \nabla^i Z_s^{\pi} ds - \int_t^T    \nabla^i Z_s^{\pi} dW_s\\
 \nabla X_t &= I_d + \int_0^t \partial_x b \nabla X_s ds.
\end{align*}
Thanks to this representation of the process $Z$, we can now apply the same strategy than in \cite{Richou-12} to show that 
 \begin{equation}
 |Z_t^{\pi}| \leqslant C(1+\sup_{s \in [0,t]} |X_s|^r) \quad \quad d\mathbb{P} \otimes dt \textrm{ a.e.}
 \end{equation}
where $C$ only depends on constants appearing in (F1)-(B2) and does not depend on $\pi$ nor on $N$. We emphasize the fact that this is possible due to the uniform (in $N$ and $\pi$) bound on $\sum_{i=0}^n |\partial_{x_i} h^{\pi}(x)|$ and $ \sum_{i=0}^n |\partial_{x_i} g^{\pi}(x,y,z)|$. When $g$ and $b$ are not smooth we can obtain the same result by a standard smooth approximation.

3. Since $h^{\pi}$ tends to $h$ and $g^{\pi}$ tends to $g$, recalling we have a Lipschitz generator, we can use a standard stability result to get that $(Y^{\pi},Z^{\pi}) \rightarrow (Y,Z)$ in $\mathcal{S}^2 \times \mathcal{M}^2$ and so 
 \begin{equation}
 |Z_t| \leqslant C(1+\sup_{s \in [0,t]} |X_s|^r) \quad \quad d\mathbb{P} \otimes dt \textrm{ a.e.}
 \end{equation}
\eproof
\begin{rem}
 \begin{itemize}
 \item  The case $r=1$ can be also tackled with extra assumptions as in \cite{Richou-12}. More precisely, we have to assume that $K_h$, $K_g$ and $T$ are small enough to ensure exponential integrability of the terminal condition and the random part of the generator. These extra assumptions are natural when we are looking for the existence of a solution, see e.g. \cite{Briand-Hu-06}.
 \item The estimate \ref{est:Z path dependent bruit additif} is interesting in itself and can be useful in many situations. For example, we can adapt the proof to obtain the same kind of estimate in a super-quadratic setting, as in \cite{Richou-12}, and then obtain an existence and uniqueness result for path-dependent super-quadratic  BSDEs. We can also use this estimate to get an explicit error bound when we consider a truncated (in $z$) approximation of the BSDE in order to deal with BSDE numerical approximation schemes (see Section 5 in  \cite{Richou-12}). See also \cite{Bender-14} for a possible application of this kind of estimate to BSDEs driven by Gaussian Processes.
\end{itemize}
\end{rem}

\begin{prop}
\label{prop:unicitepathdep1}
 We assume that Assumptions (F1) and (B2) hold. Then the BSDE \eqref{BSDE path dep 1} admits a unique solution $(Y,Z)$ satisfying
\begin{equation*}
 e^{K_z|Y^*|} \in \bigcap_{p >1} L^p.
\end{equation*}
\end{prop}
\proof
We start by considering the BSDE
\begin{equation}
\label{eq:EDSRY1Z1}
Y_t^1 = h((X_s)_{s \in [0,T]}) + \int_t^T g((X_{u \wedge s})_{u \in [0,T]},Y_s^1,Z_s^1) ds -\int_t^T Z_s^1 dW_s.
\end{equation}
Using Proposition \ref{prop:majorationZ}, we have the existence of a solution $(Y^1,Z^1) \in \mathcal{S}^2 \times \mathcal{M}^2$ to equation \eqref{eq:EDSRY1Z1} such that 
\begin{equation}
\label{est:Z1}
|Z_t^1| \leqslant C(1+\sup_{s \in [0,t]} |X_s|^r), \quad d\mathbb{P} \otimes dt \textrm{ a.e.}
\end{equation}
Now we introduce a new BSDE
\begin{eqnarray*}
\nonumber
 Y_t^2 &=& \xi+\int_t^T f(s,Y_s^1+Y_s^2,Z_s^1+Z_s^2)ds - \int_t^T Z_s^2 dW_s^{\mathbb{Q}}\\
 &&+\int_t^T \left\{g((X_{u \wedge s})_{u \in [0,T]},Y_s^1+Y_s^2,Z_s^1+Z_s^2)- g((X_{u \wedge s})_{u \in [0,T]},Y_s^1,Z_s^1)\right\}ds\\
 &&-\int_t^T Z_s^2 \nabla_z g((X_{u \wedge s})_{u \in [0,T]},Y_s^1,Z_s^1) ds,
\end{eqnarray*}
where $dW^{\mathbb{Q}}_s = dW_s-\nabla_z g((X_{u \wedge s})_{u \in [0,T]},Y_s^1,Z_s^1) ds$. By using Novikov's condition, { there} exists a probability $\mathbb{Q}$ under which $W^{\mathbb{Q}}$ is a Brownian motion.  Then, \cite{Kobylanski-00} gives us the existence of a solution $(Y^2,Z^2) \in \mathcal{S}^{\infty}(\mathbb{Q}) \times \mathcal{M}^2(\mathbb{Q})$ to the previous BSDE such that $Z^2 \in BMO(\mathbb{Q})$. Now we can remark that $(Y,Z):=(Y^1+Y^2,Z^1+Z^2)$ is a solution of \eqref{BSDE path dep 1}. We denote 
$$F(t,y,z):= f(t,y,z)+g((X_{u \wedge t})_{u \in [0,T]},y,z)$$
and we get
\begin{eqnarray*}
 \mathcal{E}_T &:=& e^{\int_0^T \nabla_z F(s,Y_s,Z_s)dW_s-\frac{1}{2} \int_0^T |\nabla_z F(s,Y_s,Z_s)|^2 ds} = e_1e_2e_3,
\end{eqnarray*}
with 
\begin{eqnarray*}
 e_1 &=& e^{\int_0^T \left( \nabla_z F(s,Y_s,Z_s)-\nabla_z F(s,Y_s,Z_s^1)\right)dW_s-\frac{1}{2} \int_0^T |\nabla_z F(s,Y_s,Z_s)-\nabla_z F(s,Y_s,Z_s^1)|^2 ds},\\ 
 e_2 &=& e^{\int_0^T \nabla_z F(s,Y_s,Z_s^1)dW_s-\frac{1}{2} \int_0^T |\nabla_z F(s,Y_s,Z_s^1)|^2 ds},\\
 e_3 &=& e^{- \int_0^T \langle \nabla_z F(s,Y_s,Z_s),\nabla_z F(s,Y_s,Z_s)-\nabla_z F(s,Y_s,Z_s^1) \rangle ds}.
\end{eqnarray*}
We will study the integrability of these terms. First of all, we can remark that
\begin{equation}
\label{est:nablazF}
\int_0^T|\nabla_z F(s,Y_s,Z_s^1)|^2ds \leqslant C(1+\sup_{s \in [0,T]} |X_s|^{2r})
\end{equation}

due to (B2)-1, (B2)-4 and \eqref{est:Z1}. By using Novikov's condition and classical estimates on exponential moments of SDEs, it implies that
\begin{equation}
 \label{est:e2danstouslesLp}
 e_2 \in \bigcap_{p\geqslant 1} L^p .
\end{equation}
For same reasons we have
\begin{equation}
 \label{eq:chgtvariabledansLP}
 e^{\int_0^T \nabla_z g((X_{u \wedge s})_{u \in [0,T]},Y_s^1,Z_s^1) dW_s-\frac{1}{2} \int_0^T |\nabla_z g((X_{u \wedge s})_{u \in [0,T]},Y_s^1,Z_s^1)|^2 ds} \in \bigcap_{p\geqslant 1} L^p.
\end{equation}
Since $|\nabla_z F(s,Y_s,Z_s)-\nabla_z F(s,Y_s,Z_s^1)| \leqslant 2K_z |Z_s^2|$, then we obtain
$$\int_0^. \left(\nabla_z F(s,Y_s,Z_s)-\nabla_z F(s,Y_s,Z_s^1)\right)dW_s \in BMO(\mathbb{Q})$$
and so there exists $\ell >1$ such that $e_1 \in L^{\ell}(\mathbb{Q})$. By using \eqref{eq:chgtvariabledansLP}
and H\"older inequality we get that $e_1 \in L^{1+\frac{\ell-1}{2}}$. We can also observe that,
$$ e_3 \leqslant e^{- \int_0^T \langle \nabla_z F(s,Y_s,Z_s^1),\nabla_z F(s,Y_s,Z_s)-\nabla_z F(s,Y_s,Z_s^1) \rangle ds}$$
and, by using Young inequality, 
\begin{eqnarray*}
  |\langle \nabla_z F(s,Y_s,Z_s^1),\nabla_z F(s,Y_s,Z_s)-\nabla_z F(s,Y_s,Z_s^1) \rangle| \leqslant \frac{1}{4\varepsilon} |\nabla_z F(s,Y_s,Z_s^1)|^2 +\varepsilon K_z^2|Z_s^2|^2,
\end{eqnarray*}
for all $\varepsilon >0$. It implies that
\begin{eqnarray*}
 e_3 &\leqslant& e^{\frac{C}{\varepsilon}(1+\sup_{s \in [0,T]} |X_s|^{2r}) +\varepsilon K_z^2\int_0^T |Z_s^2|^2ds}.
\end{eqnarray*}
Since $\int_0^. Z_s^2 dW_s$ is BMO we can apply the John-Nirenberg inequality (see \cite{Kazamaki-94}) and by using Cauchy-Schwarz inequality and classical estimates on exponential moments of SDEs we get
\begin{equation}
 \label{est:e3danstouslesLp}
 e_3 \in \bigcap_{p>1} L^p .
\end{equation}
Finally, by using \eqref{est:e2danstouslesLp}, \eqref{est:e3danstouslesLp} and the estimate $e_1 \in \cup_{p>1} L^p$, we get that 
$$\mathcal{E}_T  \in \bigcup_{p>1} L^p.$$
We just have to apply Theorem \ref{resultat unicite 3} to conclude.
\eproof

\subsection{A second unbounded case}

In this subsection we consider a more general SDE
\begin{equation}
  \label{SDE 2}
  X_t = x + \int_0^t b(X_s)ds+\int_0^t \sigma(X_s) dW_s, \quad 0 \leqslant t \leqslant T,
\end{equation}
where $b$ and $\sigma$ satisfies classical assumptions:
\paragraph{(F2)}
$ $
\begin{enumerate}
 \item $b: \mathbb{R}^d \rightarrow \mathbb{R}^d$ and $\sigma: \mathbb{R}^d \rightarrow \mathbb{R}^{d \times d}$ are Lipschitz functions: for all $(x,x') \in \mathbb{R}^d \times \mathbb{R}^d$ we have $|b(x)-b(x')| \leqslant K_b |x-x'|$ and $|\sigma(x)-\sigma(x')|\leqslant K_{\sigma} |x-x'|$.
 \item $\sigma$ is bounded by $|\sigma|_{\infty}$.
\end{enumerate}

Now we want to study the same BSDE \eqref{BSDE path dep 1}
under following assumptions:
\paragraph{(B3)}
$ $ 
\begin{enumerate}
 \item $$|\xi|_{L^{\infty}}+\left|\int_0^T |f(s,0,0)|ds\right|_{L^{\infty}}+ \left| \sup_{s \in[0,T],y \in \mathbb{R}, \mathbf{x} \in C([0,T],\mathbb{R}^d)} |\nabla_z f (s,y,0)|+ |\nabla_z g(\mathbf{x},y,0)|  \right|_{L^{\infty}} <+\infty$$
 \item there exists $C>0$ such that, for all $s\in [0,T]$, $y \in \mathbb{R}$ and $z \in \mathbb{R}^{1\times d}$,
 $$|f(t,y,z)| \leqslant C.$$
 \item There exist $K_h>0$ and $K_g>0$ such that, for all $\mathbf{x}, \mathbf{\tilde{x}} \in C([0,T],\mathbb{R}^d)$, $y \in \mathbb{R}$, $z \in \mathbb{R}^{1 \times d}$,
 $$|h(\mathbf{x}) -h({\mathbf{\tilde{x}}})| \leqslant K_h|\mathbf{x} - \mathbf{\tilde{x}}|_{\infty},$$
 $$|g(\mathbf{x},y,z) -g({\mathbf{\tilde{x}},y,z})| \leqslant K_g|\mathbf{x} - \mathbf{\tilde{x}}|_{\infty},$$
 \item (B1) holds true for $(s,y,z) \mapsto f(s,y,z)$ and $(s,y,z) \mapsto g((X_{u \wedge s})_{u \in [0,T]},y,z)$.
\end{enumerate}

Firstly we give a general lemma.
\begin{lem}
\label{lem:ZBMOZborne}
We assume that (B1) is in force and 
\begin{itemize}
 \item $\mathbb{E}\left[ |\xi|^2 + \int_0^T |f(t,0,0)|^2dt\right] <+\infty$,
 \item $\xi \in \mathbb{D}^{1,2}$ and for all $(y,z) \in \mathbb{R} \times \mathbb{R}^{1 \times d}$, $f(.,y,z) \in \mathbb{L}_{1,2}$,
 \item $\left|\sup_{y \in \mathbb{R},s \in [0,T]} |\nabla_z f(s,y,0)| \right|_{L^{\infty}}\leqslant M_z.$
\end{itemize}
Let us consider an auxiliary BSDE
\begin{eqnarray}
 \label{eq:EDSR-RS}
 R_t &=&\xi +\int_t^T f(s,R_s,0)ds- \int_t^T S_s dW_s,
\end{eqnarray}
with a unique solution $(R,S) \in \mathcal{S}^2 \times \mathcal{M}^2$. 
If 
\begin{equation}
\label{hyp:DxiDfbornes}
 \left| \sup_{t \in [0,T]} \mathbb{E}_t \left[ |D_t \xi| +\int_t^T |D_t f(s,y,0)|_{y=R_s}ds\right]\right|_{L^{\infty}} <+\infty,
\end{equation}
then $S$ is $d\mathbb{P} \otimes dt$ a.e. bounded and there exists a solution $(Y,Z) \in \mathcal{S}^2 \times \mathcal{M}^2$ of \eqref{BSDE} such that  $\int_0^. Z_s dW_s$ is BMO.

If moreover we have, for all $p >1$,
\begin{equation*}
 \left| \sup_{t \in [0,T]} \mathbb{E}_t \left[  |D_t \xi|^p + \left(\int_t^T |D_t f(s,y,z)|_{y=Y_s,z=Z_s}ds\right)^p \right]\right| <+\infty,
\end{equation*}
then there exists a solution $(Y,Z) \in \mathcal{S}^2 \times \mathcal{M}^2$ of \eqref{BSDE} such that $Z$ is $d\mathbb{P} \otimes dt$ a.e. bounded.
\end{lem}
\proof
Let us assume that $f$ and $\xi$ satisfy assumptions of Proposition 5.3 in \cite{ElKaroui-Peng-Quenez-97} (smoothness and integrability assumptions). Then we can differentiate (in the Malliavin sense) BSDE \eqref{eq:EDSR-RS}: We obtain, for all $t \in [0,T]$,
\begin{equation*}
 D_t R_t = \mathbb{E}_t \left[ e^{\int_t^T \nabla_y f(s,R_s,0)ds} D_t \xi + \int_t^T e^{\int_t^s \nabla_y f(r,R_r,0)dr} (D_t f(s,y,0))_{y=R_s}ds\right]
\end{equation*}
and a version of $S$ is given by $(D_t R_t)_{t \in [0,T]}$.
Thus we get that there exists $C>0$ such that, for all $t \in [0,T]$,
\begin{equation}
 \label{ineg:borne sur S}
\abs{S_t} = \abs{D_t R_t}\leqslant e^{K_y T} \mathbb{E}_t \left[|D_t \xi|+\int_t^T |D_t f(s,y,0)|_{y=R_s} ds\right] \leqslant C.
\end{equation}
When $f$ is not smooth enough and $f$, $\xi$ are not enough integrable, we can show by a standard approximation procedure that inequality \eqref{ineg:borne sur S} stays true $d\mathbb{P} \otimes dt$ a.e.  
Now we consider the following BSDE:
\begin{equation}
\label{eq:EDSR-UV}
U_t = \int_t^T f(s,U_s+R_s,V_s+S_s)-f(s,R_s,0) ds - \int_t^T V_s dW_s.
\end{equation}
If we set $\Psi(s,u,v) := f(s,u+R_s,v+S_s)-f(s,R_s,0)$, then, by using \eqref{ineg:borne sur S} and assumptions of the Lemma on $f$, we have, for all $s \in [0,T]$, $u,u' \in \mathbb{R}$, $v,v' \in \mathbb{R}^{1 \times d}$, 
\begin{itemize}
 \item $|\Psi(s,u,v)-\Psi(s,u',v)| \leqslant |v-v'|$, 
 \item $|\Psi(s,u,v)-\Psi(s,u,v')| \leqslant \left( M_z + K_z(|z|+|z'|+2C)\right)|z-z'|$,
 \item $|\Psi(s,u,v)| \leqslant  K_y |u| + (M_z+K_z|v+S_s|) |v+S_s| \leqslant C(1+|u|+|v|^2).$
\end{itemize}
By applying results of \cite{Kobylanski-00} we obtain a unique solution $(U,V) \in \mathcal{S}^{\infty}\times \mathcal{M}^2$ and moreover $\int_0^. V_s dW_s$ is BMO. Finally,  we can remark that $(Y,Z):= (U+R,V+S)$ is a solution of BSDE \eqref{BSDE}. So, since $S$ is bounded, $\int_0^. Z_s dW_s$ is BMO.

Concerning the boundedness of $Z$, we just have to adapt the proof of Theorem 3.6 in \cite{Richou-12} in a non Markovian framework which does not create { any difficulty}. For the reader convenience, we only sketch the proof and we refer to \cite{Richou-12} for further details. We start by assuming that $f$ and $\xi$ satisfy assumptions of Proposition 5.3 in \cite{ElKaroui-Peng-Quenez-97} (smoothness and integrability assumptions). Then we can differentiate (in the Malliavin sense) BSDE \eqref{BSDE}: We obtain, for all $t \in [0,T]$, $u \in [0,T]$,
\begin{eqnarray*}
D_u Y_t &=& D_u \xi + \int_t^T (D_u f(s,y,z))_{y=Y_s,z=Z_s} + D_u Y_s \nabla_y f(s,Y_s,Z_s) +D_u Z_s\nabla_z f(s,Y_s,Z_s) ds\\
&&- \int_t^T D_u Z_s dW_s\\
&=& D_u \xi e^{\int_t^T \nabla_y f(s,Y_s,Z_s) ds} +\int_t^T  e^{\int_s^T \nabla_y f(r,Y_r,Z_r) dr} (D_u f(s,y,z))_{y=Y_s,z=Z_s} ds - \int_t^T D_u Z_s dW_s^{\mathbb{Q}}, 
\end{eqnarray*}
where $dW_s^{\mathbb{Q}} : = dW_s - \nabla_z f(s,Y_s,Z_s) ds$ and a version of $Z$ is given by $(D_t Y_t)_{t \in [0,T]}$. Thanks to assumptions on the growth of $\nabla_z f$ and the fact that $\int_0^. Z_s dW_s$ is BMO, we know that there exists a probability $\mathbb{Q}$ under which $W^{\mathbb{Q}}$ is a Brownian motion. It implies the following estimate
\begin{equation*}
 |Z_t| = |D_t Y_t| \leqslant e^{K_y T} \mathbb{E}^{\mathbb{Q}}_t \left[ \left|D_t \xi\right| + \int_t^T |D_t f(s,y,z)|_{y=Y_s,z=Z_s} ds\right].
\end{equation*}
Then, we use once again BMO properties of $\int_0^. Z_s dW_s$: thanks to the reverse H\"older inequality (see Kazamaki \cite{Kazamaki-94}), we can apply H\"older inequality to the previous estimate to obtain the existence of $C>0$ and $p>1$ (that depend only on constants appearing in assumptions of the Lemma) such that
\begin{equation}
\label{ineg:borne sur Z}
 |Z_t| = |D_t Y_t| \leqslant C e^{K_y T} \mathbb{E}_t \left[ \left|D_t \xi\right|^p + \left(\int_t^T |D_t f(s,y,z)|_{y=Y_s,z=Z_s} ds\right)^p\right]^{1/p}.
\end{equation}
We use \eqref{hyp:DxiDfbornes} to conclude. When $f$ is not smooth enough and $f$, $\xi$ are not enough integrable, we can show by a standard approximation procedure that inequality \eqref{ineg:borne sur Z} stays true $d\mathbb{P} \otimes dt$ a.e.  

\eproof

\begin{cor}
\label{cor:Z borne}
We consider the path-dependent framework and so we assume that $\xi=0$ and $f=0$. We also assume that Assumptions (F2)-(B3) hold true. Then there exists a solution $(Y,Z)$ of the path-dependent BSDE \eqref{BSDE path dep 1} in $\mathcal{S}^2 \times \mathcal{M}^2$ such that,
 \begin{equation*}
 |Z_t| \leqslant C, \quad d\mathbb{P} \otimes dt \,\, \text{a.e.}
 \end{equation*}
\end{cor}
\proof 
1. Let us start by the Markovian framework. Without lost of generality we can assume that $X^1_t=t$ for all $t \in [0,T]$. We assume that, for all $\mathbf{x} \in C([0,T],\mathbb{R}^d)$, $y\in \mathbb{R}$, $z \in \mathbb{R}^{1 \times d}$, we have $h(\mathbf{x})=\tilde{h}(\mathbf{x}_T)$ and $g(\mathbf{x},y,z) = \tilde{g}(\mathbf{x}_{\sup_{t \in [0,T]} \mathbf{x}^1_t},y,z)$ with $\tilde{h}:\mathbb{R^d} \rightarrow \mathbb{R}$ a $K_h$-Lipschitz function and $\tilde{g} : \mathbb{R}^d \times \mathbb{R} \times \mathbb{R}^{1 \times d} \rightarrow \mathbb{R}$ a $K_g$-Lipschitz function with respect to the first variable (uniformly in $y$ and $z$). If $\tilde{h}$, $\tilde{g}$, $b$ and $\sigma$ are smooth enough then $\tilde{h}(X_T)$ and $\tilde{g}(X_s,y,z)$ are Malliavin differentiable and the chain rule gives us
$$D_t \tilde{h}(X_T) = \nabla_x \tilde{h}(X_T) \nabla X_T (\nabla X_t)^{-1} \sigma (X_t), \quad D_t \tilde{g}(X_s,y,z) = \nabla_x \tilde{g}(X_s,y,z)\nabla X_s (\nabla X_t)^{-1} \sigma(X_t) \mathbbm{1}_{t \leqslant s}.$$
So we get, for all $p\geqslant 1$,
\begin{align*}
\mathbb{E}_t\left[|D_t \tilde{h}(X_T)|^p + (\int_t^T |D_t \tilde{g} (X_s,y,z)|ds)^p \right] &\leqslant (K_h^p   + K_g^p T^p )|\sigma|_{\infty}^p\mathbb{E}_t\left[\sup_{s \in [t,T]}|\nabla X_s (\nabla X_t)^{-1}|^p \right] \\
 &\leqslant (K_h^p   + K_g^p T^p )|\sigma|_{\infty}^p C_p
\end{align*}
with $C_p$ that only depends on $p$, $T$, $K_b$ and $K_{\sigma}$, thanks to classical estimates on SDEs. Then we just have to apply the Lemma \ref{lem:ZBMOZborne} to obtain that $Z$ is bounded with a bound that only depends on constants appearing in assumptions. When $\tilde{h}$, $\tilde{g}$, $b$ and $\sigma$ are not smooth enough we can show that this result stays true by a standard approximation procedure. 

2. To deal with the general path-dependent framework we just have to apply the same strategy than in Proposition \ref{prop:majorationZ}, we firstly consider the discrete path-dependent case and then we pass to the limit. We refer to this proof for further details.
\eproof

\begin{rem}
Corollary \ref{cor:Z borne} answers an open question in the Section 3 of \cite{Richou-12}. In light to this result a new question arise:  what happens when $g$ and $h$ are only locally Lipschitz? More precisely, does Proposition \ref{prop:majorationZ} stays true when we replace assumption (F1) by assumption (F2)? Let us remark that the answer is not clear even when $f(z) = \frac{|z|^2}{2}$, see \cite{Richou-12}.
\end{rem}

\begin{prop}
\label{prop:unicitepathdep2}
 We assume that Assumptions (F2)-(B3) hold true. Then the BSDE \eqref{BSDE path dep 1} admits a unique solution $(Y,Z)$ satisfying
\begin{equation*}
 e^{K_z|Y^*|} \in \bigcap_{p >1} L^p.
\end{equation*}
\end{prop}
\proof
The proof follows the same lines than the proof of Proposition \ref{prop:unicitepathdep1}.
\eproof


\bibliographystyle{plain}

\begin{thebibliography}{10}

\bibitem{Bender-14}
C.~Bender.
\newblock Backward {SDE}s driven by {G}aussian processes.
\newblock {\em Stochastic Process. Appl.}, 124(9):2892--2916, 2014.

\bibitem{Briand-Elie-13}
P.~Briand and R.~Elie.
\newblock A simple constructive approach to quadratic {BSDE}s with or without
  delay.
\newblock {\em Stochastic Process. Appl.}, 123(8):2921--2939, 2013.

\bibitem{Briand-Hu-06}
P.~Briand and Y.~Hu.
\newblock B{SDE} with quadratic growth and unbounded terminal value.
\newblock {\em Probab. Theory Related Fields}, {\bf 136}(4):604--618, 2006.

\bibitem{Briand-Hu-08}
P.~Briand and Y.~Hu.
\newblock Quadratic {BSDE}s with convex generators and unbounded terminal
  conditions.
\newblock {\em Probab. Theory Related Fields}, {\bf 141}(3-4):543--567, 2008.

\bibitem{Delbaen-Hu-Richou-13}
F.~Delbaen, Y.~Hu, and A.~Richou.
\newblock On the uniqueness of solutions to quadratic {BSDE}s with convex
  generators and unbounded terminal conditions: the critical case.
\newblock arXiv:1303.4859v1, to appear in {D}iscrete and {C}ontinuous
  {D}ynamical {S}ystems. {S}eries {A}.

\bibitem{Delbaen-Hu-Richou-09}
F.~Delbaen, Y.~Hu, and A.~Richou.
\newblock On the uniqueness of solutions to quadratic {BSDE}s with convex
  generators and unbounded terminal conditions.
\newblock {\em Ann. Inst. Henri Poincar\'e Probab. Stat.}, 47(2):559--574,
  2011.

\bibitem{ElKaroui-Peng-Quenez-97}
N.~El~Karoui, S.~Peng, and M.~C. Quenez.
\newblock Backward stochastic differential equations in finance.
\newblock {\em Math. Finance}, {\bf 7}(1):1--71, 1997.

\bibitem{Hu-Ma-04}
Y.~Hu and J.~Ma.
\newblock Nonlinear {F}eynman-{K}ac formula and discrete-functional-type
  {BSDE}s with continuous coefficients.
\newblock {\em Stochastic Process. Appl.}, 112(1):23--51, 2004.

\bibitem{Kazamaki-94}
N.~Kazamaki.
\newblock {\em Continuous exponential martingales and {BMO}}, volume~{\bf 1579}
  of {\em Lecture Notes in Mathematics}.
\newblock Springer-Verlag, Berlin, 1994.

\bibitem{Kobylanski-00}
M.~Kobylanski.
\newblock Backward stochastic differential equations and partial differential
  equations with quadratic growth.
\newblock {\em Ann. Probab.}, {\bf 28}(2):558--602, 2000.

\bibitem{Richou-12}
A.~Richou.
\newblock Markovian quadratic and superquadratic {BSDE}s with an unbounded
  terminal condition.
\newblock {\em Stochastic Process. Appl.}, 122(9):3173 -- 3208, 2012.

\bibitem{Zhang-04}
J.~Zhang.
\newblock A numerical scheme for {BSDE}s.
\newblock {\em Ann. Appl. Probab.}, {\bf 14}(1):459--488, 2004.

\end{thebibliography}
\def\cprime{$'$}

\end{document}